\documentclass[12pt]{article}

\usepackage{amsmath,amssymb,amsthm,mathrsfs, dcpic, pictexwd}
\usepackage[margin=1in]{geometry}

\newtheorem{theorem}{Theorem}[section]

\newtheorem{prop}[theorem]{Proposition}

\newtheorem{example}[theorem]{Example}

\newtheorem{remark}{Remark}

\numberwithin{equation}{section}

\def\PP{{\mathbb{P}}}

\def\fff{{\mathbb{F}}}

\def\nnn{\mathbb{N}}
\def\Sz{{\rm{Sz}}}
\def\Div{{\rm{Div}}}

\newcommand{\oper}[1]{\operatorname{#1}} 
\newcommand{\divi}{\operatorname{Div}}
\newcommand{\supi}{\operatorname{Supp}}

\newcommand{\fqq}{\mathbb{F}_q}
\newcommand{\fqqq}{\mathbb{F}_{q^4}}

\newcommand{\divv}{\operatorname{div}}

\newcommand{\fe}{\mathbb{F}_{q^4}}
\newcommand{\cod}[2]{C_{\mathcal{L}}(#1,#2)}
\newcommand{\valu}[2]{\nu_{#1}(#2)}
\newcommand{\res}[2]{\operatorname{res}_{#1}(#2)}

\newcommand{\pif}{P_\infty}

\newcommand{\BB}{^{2}B_{2}}

\newcommand{\fgq}{\mathbb{F}_{q}}

\newcommand{\fg}[1]{\mathbb{F}_{#1}}
\newcommand{\fgqc}{\bar{\mathbb{F}}_2}

\newcommand{\sfrac}[2]{#1/#2}

\newcommand{\placeextension}{|}

\begin{document}

\title{Suzuki-invariant codes from the Suzuki curve}

\author{Abdulla Eid, Hilaf Hasson, Amy Ksir, Justin Peachey\thanks{This work was conducted at the Mathematics Department, Colorado State University, Summer 2011, funded by NSF grant DMS-11-01712}}


\date{\today}

\maketitle

\begin{abstract}
In this paper we consider the Suzuki curve $y^q + y = x^{q_0}(x^q + x)$ over the field with $q = 2^{2m+1}$ elements.  The automorphism group of this curve is known to be the Suzuki group $\Sz(q)$ with $q^2(q-1)(q^2+1)$ elements.  We construct AG codes over $\fqqq$ from a $\Sz(q)$-invariant divisor $D$, giving an explicit basis for the Riemann-Roch space $L(\ell D)$ for $0 < \ell \leq q^2-1$. These codes then have the full Suzuki group $\Sz(q)$ as their automorphism group.  These families of codes have very good parameters and are explicitly constructed with information rate close to one. The dual codes of these families are of the same kind if $2g-1 \leq \ell \leq q^2-1$.
\end{abstract}


\section{Introduction}
The Suzuki curve has been a source of very good error-correcting codes. Codes constructed from the Suzuki curve have been studied, for example, in \cite{CD}, \cite{DPark}, \cite{HS} (one-point codes), and \cite{M} (two-point codes) and shown to have very good parameters. Furthermore, the Suzuki curve has a very large automorphism group for its genus, namely the Suzuki group $\Sz(q)$ = $\BB$ of order $q^2(q-1)(q^2+1)$. The one-point and two-point codes previously studied had automorphism groups which were not the full Suzuki group.  In this paper, we construct a family of codes on the Suzuki curve with the full Suzuki group as its group of automorphisms.  We find that our codes also have very good parameters.

The outline of our paper is as follows: In Section~\ref{prelim} we start with some preliminaries about the Suzuki curve. In Section~\ref{RiemannRoch} we give an explicit basis for the Riemann-Roch space $L(\ell D)$ for $0 < \ell \leq q^2-1$, where the divisor $D$ is the sum of all $\mathbb{F}_q$--rational points of the Suzuki curve. In Section~\ref{AGCodes} we construct families of AG--codes with good parameters, and with the full Suzuki group as automorphism group. These families are significant because they are explicitly constructed in a polynomial--time with rate close to one. In Section~\ref{Dual} we find the dual codes of the codes constructed in Section~\ref{AGCodes} and we find the conditions of these codes to be of the same kind, isodual, and iso--orthogonal.

The authors would like to thank Rachel Pries for organizing the workshop on rational points on Suzuki Curves in which this paper was conceived.
\section {Preliminaries}\label{prelim}
Let $m \ge 1$ be an integer, $q_0:=2^m$, $q:=2^{2m+1}=2q_0^2$, and let $X_m$ denote the smooth projective curve with affine plane equation
\begin{equation}\label{eq2:fm}
 y^q+y=x^{q_0}(x^q+x)
\end{equation}
over $\fgq$. Then, $X_m$ has a singular projective plane model $Y_m$ in $\mathbb{P}^{2}_{\fgqc}$ with the homogeneous equation
\[
 y^qt^{q_0}+yt^{q+q_0-1}=x^{q+q_0}+x^{q_0+1}t^{q-1}
\]
in homogeneous coordinates $[t:x:y]$. This curve has been studied, for example, in \cite{DL} and \cite{Henn} and it has been shown in \cite{TGK} that the curve has a smooth projective embedding in $\mathbb{P}^4$.  Moreover, $X_m$ has a very large automorphism group for its genus, namely the Suzuki group $\Sz(q)$ of order $q^2(q-1)(q^2+1)$. As such, $X_m$ is known as the Suzuki curve.  We summarize these properties as well as several others shown in \cite{FTor1} and \cite{HS} in the following proposition:
\begin{prop}\label{prop2:1}
Let $m \ge 1$ be an integer, $q_0:=2^m$, $q:=2^{2m+1}=2q_0^2$, and let $X_m$ denote the Suzuki curve.  Then,
\begin{enumerate}
  \item
  The smooth projective curve $X_m$ has a single point $P_{\infty}$ above the singularity at infinity $[0:0:1]$ of $Y_m$.
  \item The genus of $X_m$ is $g:=q_0(q-1)$.
  \item The number of $\fgq$-rational points is $q^2+1$, which is maximal as shown by the Serre bound.
  \item The Suzuki curve $X_m$ is the unique curve (up to $\fgq$-isomorphism) with properties (2) and (3) above.
  \item The automorphism group of $X_m$, as well as of $X_m\times_{\mathbb{F}_q}\bar{\mathbb{F}}_2$, is the Suzuki group $\Sz(q)$ = $\BB$ of order $q^2(q-1)(q^2+1)$.
  \item The functions $x$, $y$, $z:=x^{2q_0+1}-y^{2q_0}$, and $w:=xy^{2q_0}-z^{2q_0}$ are regular outside $P_\infty$ with pole orders at $P_\infty$ given by $q,q+q_0,q+2q_0,$ and $q+2q_0+1$ respectively.
  \item The functions $t$, $x$, $y$, $z$ and $w$ give a smooth embedding of $X_m$ into $\PP^4$. 
 \end{enumerate}
\end{prop}
The number $N_j(X_m)$ of $\fff_{q^j}$-rational points on the curve can be determined using the zeta function of the curve, or more specifically using the $L$ polynomial, which is the numerator of the zeta function, as follows. By \cite[Corollary 5.1.16]{Sti}, if the $L$-polynomial is $L(X_m,t)=\prod_{k=1}^{2g}(1-\alpha_kt)$, then
\begin{equation}
N_j(X_m)=q^j+1-\sum_{k=1}^{2g}\alpha_k^j \label{NPoint}.
\end{equation}
For the Suzuki curve, it was shown in \cite{H} that \[L(X_m, t)=(1+2q_0t+qt^2)^g.\] The roots of the polynomial $L(X_m, t)$ are $\underbrace{\alpha,\alpha,\dots,\alpha}_{g \text{ times}}$ and $\underbrace{\beta,\beta,\dots,\beta}_{g \text{ times}}$, where \[\alpha:=q_0(-1+i)\] and \[\beta:=\bar{\alpha}=q_0(-1-i).\]

Therefore, $$N_j(X_m)=q^j+1-gq_0(-1+i)-gq_0(-1-i)=q^j+1+2gq_0.$$

In the case $j=1$, we see that $N_1(X_m) = q + 1 - q(q-1) = q^2 + 1$.  
We will use these points to construct our codes.

\section{The Riemann-Roch space $\mathcal{L}(\ell D)$}
\label{RiemannRoch}

In order to construct an AG code whose automorphism group is the full automorphism group $\Sz(q)$ of $X_m$, we need to choose a divisor that is invariant under the action of $\Sz(q)$ on $X_m$.  Suzuki originally  constructed $\Sz(q)$ as a doubly transitive group acting on the curve \cite{Sz}.  So the only way to choose an invariant divisor is to take the set of \emph{all} $\fff_{q^j}$ points for some $j$.  The smallest such set of points is the set of $\fff_q$-points.

Consider the divisor $D\in \Div(X_m)$ given by the sum of all $\fff_q$-rational points of $X_m$.  These are the points $P_{\alpha, \beta}$ with affine coordinates $(\alpha,\beta)$ for any $\alpha$ and $\beta$ in $\fff_q$, plus the point at infinity.  Thus
\[D=P_{\infty}+\sum_{\alpha, \beta \in \fff_q} P_{\alpha,\beta}.\] 
Since there are $q^2 + 1$ many $\fgq$-rational points of $X_m$, $\deg(D) = q^2 + 1$. Moreover, the divisor $D$ is fixed by $\Sz(q)$, so codes based at $D$ will have $\Sz(q)$ as their automorphism group.  In this section, we prove the following theorem, finding an explicit $\fff_q$-basis for the space $\mathcal{L}(\ell D)$, where $\ell \leq q^2-1$.

\begin{theorem}\label{thm2:3}
Let $\ell \in \nnn$, $\ell \leq q^2 - 1$, and $D$ be defined to be the sum of all $\fgq$-rational points of $X_m$. Then,
\begin{equation}\label{eq:Basis}
S := \left\{\frac{x^a y^b z^c w^d}{(x^q + x)^r}:
										\begin{array}{l}
											{\small aq + b(q+q_0) + c(q + 2q_0) + d(q +2q_0 +1) \leq rq^2 + \ell,}\\
											0 \leq a \leq q - 1, 0 \leq b \leq 1, 0 \leq c \leq q_0 - 1, \\
											0 \leq d \leq q_0 - 1, 0 \leq r \leq \ell\\
										\end{array}
										\right\}
\end{equation}
is a basis for $\mathcal{L}(\ell D)$.
\end{theorem}

Note that the function $x^q + x$ vanishes at every affine point of $X_m$ and has a pole of order $q^2$ at $P_{\infty}$.  Therefore
\[
\divv(x^q+x) = -q^2P_{\infty}+\sum_{\alpha, \beta \in \fff_q} P_{\alpha,\beta}.
\]
Hence, $\ell D=\ell (q^2+1)P_\infty+\divv ((x^q+x)^\ell)$, i.e., $\ell D \sim \ell(q^2+1)P_{\infty}$. Thus, we have that $\mathcal{L}(\ell D) \simeq \mathcal{L}(\ell(q^2+1)P_\infty)$ where the $\fff_q$-isomorphism is given by  $f  \mapsto f(x^q+x)^\ell$ for $f \in \mathcal{L}(\ell D).$  Thus Theorem \ref{thm2:3} is equivalent (via $r_{\oper{new}} = \ell - r_{\oper{old}}$) to the following.

\begin{theorem}\label{thmdoodad}
Let $\ell \in \nnn$, $\ell \leq q^2 - 1$.  Then
\begin{equation} \label{eq:BasisInfinity}
S' := \left\{x^a y^b z^c w^d (x^q + x)^{r}:
										\begin{array}{l}
											{\small aq + b(q+q_0) + c(q + 2q_0) + d(q +2q_0 +1) + r q^2 \leq \ell(q^2 + 1)}\\
											0 \leq a \leq q - 1, 0 \leq b \leq 1, 0 \leq c \leq q_0 - 1, \\
											0 \leq d \leq q_0 - 1, 0 \leq r \leq \ell\\
										\end{array}
										\right\}
\end{equation}
is a basis for $\mathcal{L}(\ell(q^2+1) P_{\infty})$.
\end{theorem}

In order to prove this theorem, we recall a result in \cite{HS}. Let $\mathcal{P}\subseteq \mathbb{Z}_{\geq 0}$ be the semigroup generated by the pole orders of the functions $x$, $y$, $z$, and $w$ defined in Proposition \ref{prop2:1}. That is,
\begin{equation}
\label{eq:semigroup}
 \mathcal{P}:=\left \langle q,q+q_0,q+2q_0,q+2q_0+1\right \rangle \subseteq \mathbb{Z}_{\geq 0}.
\end{equation}
Proposition 1.6 in \cite{HS} is equivalent to the following:
\begin{prop}\label{propeptyprop}(\cite{HS})
For every integer $j$, $$\oper{dim}_{\mathbb{F}_q}(\mathcal{L}(jP_{\infty}))=\#\{n\in\mathcal{P}|n\leq j\}.$$
\end{prop}
We are now ready for the proof.
\begin{proof}(Theorem \ref{thmdoodad})
Let $f = x^a y^b z^c w^d (x^q+x)^r$ be an element of $S'$, and let $v_{\infty}$ be the discrete valuation corresponding to the point $P_{\infty}$.  Then
$$v_{\infty} (f) = -\left[aq + b(q+q_0) + c(q + 2q_0) + d(q + 2q_0 + 1) + r q^2 \right],$$
and $f$ has no other poles.
Thus the first inequality in the definition of $S'$ shows that $S' \subseteq \mathcal{L}(\ell(q^2+1) P_{\infty})$.
Thus, in light of Proposition \ref{propeptyprop}, it suffices to show that for every $n \in \mathcal{P}$ such that  $n \leq \ell(q^2+1)$,
$S'$ contains exactly one function with a pole of order $n$ at $P_{\infty}$.
First we show that  the valuations at $P_\infty$ of the functions in $S'$ are distinct.
Suppose that  $F_1$ and  $F_2$ in $S'$ had the same valuation at infinity, where $F_1 = x^{a_1} y^{b_1} z^{c_1} w^{d_1}(x^q + x)^{r_1}$ and $F_2 = x^{a_2} y^{b_2} z^{c_2} w^{d_2}(x^q + x)^{r_2}$. Then
\begin{align}\label{eq:infvals}
a_1q + b_1(q+q_0) + c_1(q + 2q_0) + d_1(q + 2q_0 + 1)+ r_1q^2 &= \nonumber \\
a_2q + b_2(q+q_0) + c_2(q + 2q_0) + d_2(q + 2q_0 + 1) + r_2q^2.
\end{align}

We consider (\ref{eq:infvals}) modulo $q_0$. Then,
\[d_1 \equiv d_2 \pmod{q_0}.\]
Since $1 \leq d_1, d_2 \leq q_0 - 1,$ it must be that $d_1 = d_2$.

Next, we consider (\ref{eq:infvals}) modulo $2q_0$. Then,
\[b_1q_0 + d_1 \equiv b_2q_0 + d_1 \pmod{2q_0}.\]
Note that $0 \leq b_1, b_2 \leq 1$. Thus, is must be that $b_1 = b_2$.

Next, consider (\ref{eq:infvals}) modulo $q$. Since $d_1=d_2$ and $b_1=b_2$, we get
\[ 2c_1q_0 \equiv  2c_2q_0 \pmod{q}\]
and therefore $c_1 \equiv c_2 \pmod{q_0}$.  Since $0 \leq c_1, c_2 \leq q_0 - 1$, it must be the case that $c_1 = c_2$.

Finally, consider (\ref{eq:infvals}) modulo $q^2$. Then, since $b_1 = b_2$, $c_1 = c_2$, $d_1 = d_2$, we have
\[a_1 q \equiv a_2 q \pmod{q^2}.\]
Note that $0 \leq a_1, a_2 \leq q-1$. Thus, it must be that $a_1 = a_2$. This also shows that $r_1 = r_2$. We conclude that if $v_\infty (F_1) = v_\infty (F_2)$, then $F_1 = F_2$.

Now we must show that if $n \leq \ell(q^2+1)$ is an element of $\mathcal{P}$, then there is a function in $S'$ with pole order $n$ at $P_{\infty}$.  Let $n$ be such an element.  By definition,
\begin{equation} \label{eq:originaln}
n= aq + b(q+q_0) + c(q+2q_0) + d(q+2q_0+1)
\end{equation}
for some positive integers $a$, $b$, $c$, $d$.
We need to show that there are $a'$,$b'$,$c'$,$d'$, and $r$ such that
\begin{equation}
\label{eq:desiredn}
n= a'q + b'(q+q_0) + c'(q+2q_0) + d'(q+2q_0+1) + rq^2
\end{equation}
and $0 \leq a' \leq q-1$, $0 \leq b' \leq 1$, $0 \leq c' \leq q_0-1$,  $0 \leq d' \leq q_0-1$, and $0 \leq r \leq \ell$.

Let $d'$ be the remainder when $n$ is divided by $q_0$.   Then $d'$ will be in the correct range.  Let
$$ n_d = \frac{n - d'(q + 2q_0 +1)}{q_0}.$$
Let $b'$ be the remainder when $n_d$ is divided by $2$.  Again, $b'$ will be in the correct range.  Let
$$ n_b = \frac{n_d - b'(2 q_0 + 1)}{2}.$$
Let $c'$ be the remainder when $n_b$ is divided by $q_0$.  Now $0 \leq c' \leq q_0-1$.  Let
$$ n_c = \frac{n_b - c'(q_0 + 1)}{q_0}.$$
Finally, let $a'$ be the remainder when $n_c$ is divided by $q$, so that $0 \leq a' \leq q-1$, and let
$$ r = \frac{n_c-a'}{q}.$$

Then we can put these back together as follows:
\begin{eqnarray*}
n &=& n_d q_0 + d'(q + 2q_0 + 1) \\
  &=& (2 n_b + b'(2q_0 +1)) q_0 + d'(q + 2q_0 + 1) \\
 &=& 2 q_0 n_b + b'(q + q_0) + d'(q + 2q_0 + q) \\
 &=& 2q_0 (q_0 n_c + c'(q_0 + 1)) +b'(q + q_0) + d'(q + 2q_0 + q) \\
 &=& n_c q + c' (q + 2q_0) +b'(q + q_0) + d'(q + 2q_0 + q) \\
&=&  r q^2 + a' q + c' (q + 2q_0) +b'(q + q_0) + d'(q + 2q_0 + q).
\end{eqnarray*}
What remains is to show that $r$ is in the correct range.  Since $n \leq \ell(q^2+1)$, this means that
\begin{eqnarray*}
n_d &\leq& \ell( 2 q q_0 + \frac{1}{q_0}) \\
n_b &\leq& \ell( q q_0 + \frac{1}{2 q_0}) \\
n_c &\leq& \ell( q + \frac{1}{q}) \\
r &\leq& \ell  + \frac{\ell}{q^2}.
\end{eqnarray*}
Since $r$ is an integer and $\ell < q^2$, this means that $r \leq \ell$.
Finally, to see that $0 \leq r$ we need first to show that $n_d,n_b,n_c$ are positive integers.  Since $n \equiv d' \pmod {q_0}$, we have that $d \equiv d' \pmod {q_0}$ and so we can write $d-d' =t_d q_0$, for some positive integer $t_d$ (note the previous assumtion asserts that $d \geq d'$, otherwise we have $d$ already in the required range and we don't need to find $n_d$). Now we have that
\begin{align*}
n_d:&=\frac{n-d'(q+2q_0+1)}{q_0}=\frac{aq+b(q+q_0)+c(q+2q_0) + (d-d')(q+2q_0+1)}{q_0}\\
      &= a(2q_0) + b (2q_0+1) +c (2q_0+2) + t_d (q+2q_0+1)\end{align*}
which is a positive integer.

Next we have that $n_d \equiv b' \pmod 2$, so we have $b+t_d \equiv b' \pmod 2$ and so again we can write it as $b+t_b-b' = 2t_b$, for some positive integer $t_b$ ( with the same assumption as before that b is not in the required range $\{0,1\}$, so $b \geq b'$). Now we have that
\begin{align*}
n_b:&= \frac{n_d - b' (2q_0+1)}{2}=\frac{a(2q_0+ c(2q_0+2) +(b+t_d-b')(2q_0+1) + t_d (q) }{2}\\
         &= a(q_0) + c (q_0+1) + t_b (2q_0+1) + t_dq_0^2
\end{align*}
Which is again a positive integer.

Next, we have that $n_b \equiv c'\pmod {q_0} $, so we have $c+t_b\equiv c' \pmod q_0$ and we write $c+t_b-c' = t_c q_0$, for some positive integer $t_c$ and we have that
\begin{align*}
n_c :&=\frac{n_b-c'(q_0+1)}{q_0}=\frac{a(q_0) + (c+t_b-c')(q_0+1)+t_bq_0+t_dq_0^2 }{q_0}\\
       &= a+t_c(q_0+1)+t_b + t_dq_0
\end{align*}
which is a positive integer.
Finally, we have $n_c \equiv a' \pmod q$ and so we have $n_c -a' $ is multiple of $q$ and thus $r$ is a positive integer.

\end{proof}

\begin{remark}\label{remarkable}\rm

The dimension of $\mathcal{L}(\ell D)$ is given by
\begin{equation} \label{eq:dimLlD}
\dim _{\mathbb{F}_q}\mathcal{L}(\ell D) = \ell(q^2+1) - q_0(q-1) +1,
\end{equation}
which we can see in two ways.  First, since $q^2+1 > 2q_0(q-1)$, we have $\deg D > 2g$ and the result follows from the Riemann-Roch theorem.  Second, in \cite[Appendix A]{HS}, it is shown that $\#(\mathbb{N} \setminus \mathcal{P}) = q_0(q-1)$, and an analysis of their proof shows that the largest number in $\mathbb{N} \setminus \mathcal{P}$ is $2q_0(q-1)-1$.  Thus $\#S$ is the number of possibilities for $a$, $b$, $c$, $d$, and $r$, minus $\#(\mathbb{N} \setminus \mathcal{P})$.

\end{remark}

Theorem \ref{thm2:3} gives us an explicit basis for $\mathcal{L}(\ell D)$, which we use to construct Suzuki-invariant codes in the next section. 

\section{Construction and properties of the code $C(E,\ell D)$}\label{AGCodes}

As above, let $D \in \divi(X_m)$ be the sum of all $\fqq$-rational points in $X_m$. By Theorem \ref{thm2:3}, for $\ell \leq q^2-1$ the Riemann-Roch space $\mathcal{L}(\ell D)$ has the $\fqq$-basis \eqref{eq:Basis}, and by Remark \ref{remarkable} $\dim_{\fqq}\mathcal{L}(\ell D)=\ell(q^2+1)-g+1=\ell(q^2+1)-q_0(q-1)+1$.

Now to construct a Suzuki-invariant geometry code, we must choose another set of points, disjoint from $D$, which is also invariant under $\Sz(q)$.  Since we have used all of the $\fff_q$ points for $D$, we must look to points over extensions of $\fff_q$.
Consider the field extension $\fqqq$ of $\fqq$. Let the divisor $E\in \divi(X_m)$ be the sum of all $\fqqq$-points minus the sum of all $\fqq$-points. Then, we have
\[
n:=\deg(E)=N_4(X_m)-N_1(X_m),
\]
where $N_4(X_m)$ is given by the formula \eqref{NPoint}, i.e.,
\begin{align*}
N_4(X_m)&=q^4+1-g(\alpha^4+\beta^4)=q^4+1+2gq^2=q^4+1+2q_0q^2(q-1).
\end{align*}
Therefore, $n=\deg(E)=q^4+1+2q_0q^2(q-1)-(q^2+1)=q^4+2q_0q^2(q-1)-q^2$.

Since $\supi(E)\cap \supi(D)=\varnothing$ and Theorem \ref{thm2:3} provides an explicit basis for $\mathcal{L}(\ell D)$, we construct an algebraic geometry code using the divisors $E,D$ as follows. Let $P_1,\dots,P_n$ be all the points in support of $E$. Define
\[
C_{m,\ell}:=C_{\mathcal{L}}(E,\ell \cdot D)=\{ (f(P_1),f(P_2),\dots,f(P_n)) \in \fqqq^{n} \,|\, f \in \mathcal{L}(\ell\cdot D) \}.
\]
Then, we have the following theorem.

\begin{theorem}\label{prop2:4}
Consider the algebraic geometry code\[C_{m,\ell}:=C_{\mathcal{L}}(E,\ell \cdot D)\] over $\fg{q^4}$, where $\ell \leq q^2-1$. Then, $C_{m,\ell}$ is an $[n,k,d]$-linear code, where
\begin{align*}
&n=\deg(E)=q^4+2q_0q^2(q-1)-q^2,\\ &k:=\dim_{\fqq}\mathcal{L}(\ell \cdot D)=\ell(q^2+1)-q_0(q-1)+1,\\
&d \geq d^*:=n-\deg(\ell \cdot D)=n-\ell(q^2+1).
\end{align*}
Moreover, this code can correct at least $t=\lfloor \sfrac{(d^*-1)}{2} \rfloor$ errors, and has $\Sz(q)$ as its automorphism group.
\end{theorem}
\begin{remark}\rm
Let $S=\{f_1,\dots,f_k \}$ be the $\fqq$-basis for $\mathcal{L}(\ell D)$ as in Theorem \ref{thm2:3}. Then, the code $C_{m,\ell}$ has generator matrix $G_{m,\ell}:=(f_j(P_i))_{1\leq j \leq k, 1\leq i \leq n}$.
\end{remark}

\begin{proof}
The parameters $n$ and $k$ were computed above; $d$, $d^*$ and $t$ come from the general theory of AG codes.  Because we chose $D$ to be an invariant divisor, the code will have the full group $\Sz(q)$ as its automorphism group.
\end{proof}

\begin{remark}\rm
It is easy to check, using equation \eqref{NPoint}, that the number of points of $X_m$
over $\fff_{q^2}$ and $\fff_{q^3}$ are
\[ N_2(X_m) = q^2+1; \qquad N_3(X_m) = q^3 + 1 - q^2(q-1) = q^2 + 1.\]
Thus if we let $E_j$ be the sum of all $\fg{q^j}$-points minus the sum of all $\fgq$-points, then $\deg(E_2)=\deg(E_3)=0$. Therefore $E_4$, which is the $E$ we used above, is the first non-trivial case. 

In fact, the Suzuki curve is a maximal curve over $\mathbb{F}_{q^4}$, meeting the Hasse-Weil bound.
\end{remark}

We now focus on the family of codes where $\ell=q^2-1$. Denote this family by $C_m$, i.e., $C_m:=C_{m,q^2-1}=C_{\mathcal{L}}(E,(q^2-1)D)$. By Theorem \ref{prop2:4}, $C_m$ is a $[n,k,d\geq d^*]$-linear code, where
\begin{align*}
n &=q^4+2q_0q^2(q-1)-q^2,\\
k &=(q^2-1)(q^2+1)-q_0(q-1)+1=q^4-q_0(q-1),\\
d^* &=n-q^4+1=2q_0q^2(q-1)-q^2+1.
\end{align*}
Using the above, $C_m$ has information rate
\begin{align*}
R_m:=\frac{k_m}{n_m}&=\frac{q^4-q_0q+q_0}{q^4+2q_0q^2(q-1)-q^2} =\frac{16q_0^8-2q_0^3+q_0}{16q_0^8+16q_0^7-4q_0^5-4q_0^4}.
\end{align*}
Thus, as $m \to \infty$, we have $R_m\to 1$. This shows that these codes are siginficant because they have very good parameter, explicitly constructed in polynomial--time, with rate asymptotically approaches one, and in many cases cannot be achieved by Reed--Solomon codes.

\begin{example}\rm The rate gets close to one very quickly. In order to show this, let us consider the following examples:
\begin{enumerate}
 \item
Let $m=1$; thus, $q=8,q_0=2$. Then, the resulting code $C_1=C_{1,63}$ is a $[5824,4082,\geq 1729]$-linear code over $\mathbb{F}_{4096}$ and can correct up to 864 errors with information rate $R_1=0.7008$. 
\item Let $m=2$; thus, $q=2^5=32, q_0=4$. Then, the resulting code $C_2=C_{2,1023}$ is a $[1051679,1048452,\geq 3104]$-linear code over $\mathbb{F}_{1048576}$ and can correct at least 1551 errors with information rate $R_2=0.996$. 
\end{enumerate}
\end{example}

\section{Dual code}\label{Dual}

As before, let $D$ be the sum of all $\fqq$-points and the divisor $E$ is the sum of all $\fe$-rational points minus all the $\fqq$-rational points. Next, we study the dual code of the code $C_{m,\ell}:=\cod{E}{\ell D}$, where $\ell \leq q^2-1$.

Recall from \cite[Proposition 2.2.10]{Sti} that the dual of an algebraic geometry code is given by $\cod{E}{\ell D}^{\perp}=\cod{E}{E-\ell D + (\eta)}$ , where $\eta$ is a Weil differential such that $\valu{P_i}{\eta}=-1$ and $\res{P_i}{\eta}=1$, for all $i=1,2,\dots,n$. 

In order to find $\eta$ we first identify the points of $\mathbb{P}^1_{\mathbb{F}_{q^4}}$ whose fiber in $X_m$ via the map induced by $x$ has an $\mathbb{F}_{q^4}$-rational point. By Proposition \ref{prop2:1}, $P_{\infty}$ is the unique point above infinity.

We therefore focus on the points that lie in the affine patch of $X_m$ isomorphic
to the open affine $t\neq 0$ of the model $Y_m$. Note that $y^q+y-\alpha^{q_0}(\alpha^q+\alpha)$ factors completely into linear terms over $\mathbb{F}_{q^4}$ for exactly $q^3+2gq$ many $\alpha$'s in $\mathbb{F}_{q^4}$, and that for the rest of the $\alpha$'s
the polynomial factors into $q/2$ many irreducible components, each of degree $2$. By Kummer's criterion \cite[Theorem 3.3.7]{Sti} applied to the equation $y^q + y - x^{q_0}(x^q + x)=0$, this implies that there are exactly $q^3+2gq$ many $\mathbb{F}_{q^4}$-rational points of $\mathbb{A}^1_{\mathbb{F}_{q^4}}$
that completely split in $X_m$, whereas the rest of $\mathbb{F}_{q^4}$-rational points of $\mathbb{A}^1_{\mathbb{F}_{q^4}}$ don't have an $\mathbb{F}_{q^4}$-rational point in their fiber.

Let $T$ be the set of $\alpha$'s in $\mathbb{F}_{q^4}$ such that $x=\alpha$ splits, and let $t:=\prod_{\alpha\in T}(x-\alpha)$ be viewed as an element of the function field $\kappa(X_m)$ of $X_m$, and let $\eta:=dt/t$. 

By the above discussion, every $\mathbb{F}_{q^4}$-rational point $P$ of $X_m$ except for $P_{\infty}$ lies above some affine point $Q_{\alpha}:=(x-\alpha)$ of $\mathbb{P}^1_{\mathbb{F}_{q^4}}$ where $\alpha\in T$.
Therefore:
\[
 v_P(t)=e(P \placeextension Q_\alpha)v_{Q_\alpha}(t)=1\cdot v_{Q_\alpha}(\prod_{\alpha \in T}(x-\alpha))=1.
\]

And so:

$$v_P(\eta)=v_P(dt/t)=-1$$

and

$$\oper{res}_P(\eta)=\oper{res}_P(1/t)=1$$

Therefore $\eta$ satisfies the conditions of
Proposition 2.2.10 in \cite{Sti}. We will now compute $(\eta)$. Note that $t$ has zeros at all the $\fe$-rational places except at $P_\infty$, that is,
\[
 (t)_0=E+D-\pif.
\]

Let $Q_{\infty}$ denote the point at infinity of $\mathbb{P}^1_{\mathbb{F}_{q^4}}$. Then:
$$v_{P_{\infty}}(t)=v_{P_{\infty}}(\prod_{\alpha\in T}(x-\alpha))=e(P_{\infty}|Q_{\infty})v_{Q_{\infty}}(\prod_{\alpha\in T}(x-\alpha))=-q|T|=-q(q^3+2gq)$$
$$=-q^4-2gq^2.$$

Therefore $(t)_{\infty}=(q^4+2gq^2)P_{\infty}$. Therefore,
\[
 (t)=(t)_0-(t)_\infty=E+D-\pif - (q^4+2gq^2)\pif.
\]

It follows that $(\eta)=(dt/t)=(2g-2)\pif-E-D+\pif + (q^4+2gq^2)\pif$. Thus, by \cite[Proposition 2.2.10]{Sti}, the dual of $C_{m,\ell}=\cod{E}{\ell D}$ is given by $ \cod{E}{G^{\perp}}$, where
\begin{align*}
 G^{\perp}&=E-\ell D +(\eta)\\
          &=E-\ell D+ (2g-2)\pif-E-D+\pif + (q^4+2gq^2)\pif\\
          &=(-1-\ell)D+(2g-2+1+q^4+2gq^2)\pif\\
          &=(-1-\ell)D+(q^2-1+2g)(q^2+1)\pif .
\end{align*}
Since $D\sim (q^2+1)\pif$,
\begin{align*}
 G^{\perp} & \sim (-1-\ell)D+(q^2+2g-1)D\\
           & \sim (q^2+2g-1-1-\ell)D.
\end{align*}

Thus, the dual code of $\cod{E}{\ell D}$ is equivalent to the code $\cod{E}{(q^2+2g-2-\ell)D}$.

Moreover, the dual code $\cod{E}{(q^2+2g-2-\ell)D}$ is also of the form $C_{m,\ell'}=\cod{E}{\ell' D}$ if $q^2+2g-2-\ell \leq q^2-1$, i.e., whenever $2g-1 \leq \ell \leq q^2-1$. (In which case $\ell'=q^2+2g-2-\ell$.)

Thus, we obtain the following result.
\begin{prop}\label{prop2:5}
If $\ell \leq q^2-1$, then the dual code of $\cod{E}{\ell D}$ is equivalent to the code $\cod{E}{(-1-\ell+q^2+2g-1)D}$.  Moreover, if $2g-1 \leq \ell$, $C_{m, \ell}^\perp$ is of the form $C_{m,\ell'}$ for
$\ell'=q^2+2g-2-\ell$.
\end{prop}

\begin{remark}\rm
The code $\cod{E}{\ell D}$ is isodual if and only if for the Weil differential above, we have
\begin{align*}
 2\ell D-E&=(\eta)\\
 2\ell D -E&=(2g-2)\pif - E -D + \pif +(q^4+2gq^2)\pif\\
 2\ell D&=(q^4+2gq^2+2g-1)\pif-D\\
 2 \ell D&\sim (q^2+2g-1)D-D\\
 2\ell&=q^2+2g-2\\
 \ell&=\frac{q^2}{2}+g-1.
\end{align*}
Hence, we have:
\begin{enumerate}
\item $\cod{E}{\ell D}$ is isodual if and only if $\ell=q^2/2+g-1$.
\item $\cod{E}{\ell D}$ is iso-orthogonal if and only if $\ell\leq q^2/2+g-1$.
\end{enumerate}
\end{remark}

\begin{example} \rm The smallest case of an isodual code in our family is the case $m=1$ and $\ell=q^2/2+g-1=8^2/2+14-1=45$. In that case the code $C_{1,45}$ is isodual.
\end{example}

\begin{remark}\rm\
Note that since the codes $\cod{E}{\ell D}$ and $\cod{E}{\ell (q^2+1)\pif}$ are equivalent (since $D\sim (q^2+1)\pif$), the code $C_{m,\ell}$ is a one-point algebraic geometry code.

\end{remark}


\medskip
\noindent Current authors information:\\ 
Abdulla Eid: Department of Mathematics at BTC, University of Bahrain, Bahrain\\
email: {\tt aeid@uob.edu.bh}\\ \\
Hilaf Hasson: Department of Mathematics, Stanford University, Palo Alto, CA 94305, USA\\
email: {\tt hilaf@stanford.edu}\\ \\
Amy Ksir: Department of Mathematics, United States Naval Academy, Annapolis, MD 21402, USA\\
email: {\tt ksir@usna.edu}\\ \\
Justin Peachey: Independent Researcher\\ 
email: {\tt jdpeachey@gmail.com}

\end{document}